\documentclass{birkart}
\usepackage{amsmath}
\usepackage{amssymb}
\usepackage{amsfonts}
\usepackage{graphicx}
\usepackage{texdraw}
\usepackage{graphpap}

\input txdtools

 \newtheorem{thm}{Theorem}[section]
 \newtheorem{cor}[thm]{Corollary}
 \newtheorem{lem}[thm]{Lemma}
 \newtheorem{prop}[thm]{Proposition}
 \theoremstyle{definition}
 \newtheorem{defn}[thm]{Definition}
 
 \newtheorem{prob}[thm]{Problem}
 \theoremstyle{remark}
 
 \newtheorem{rem}[thm]{Remark}

\numberwithin{equation}{section}
\numberwithin{figure}{section}

\newcommand{\proj}{{\mathbf P}}
\newcommand{\poiss}{\mathbf{U}}

\newcommand{\e}{\mathrm e}

\newcommand{\C}{{\mathbb C}}
\newcommand{\D}{{\mathbb D}}

\newcommand{\T}{{\mathbb T}}
\newcommand{\R}{{\mathbb R}}
\newcommand{\Sph}{{\mathbb S}}
\newcommand{\Z}{{\mathbb Z}}

\newcommand{\EP}{m^2}

\newcommand{\ci}{\operatorname{\mathrm{ci}}}
\newcommand{\si}{\operatorname{\mathrm{si}}}

\newcommand{\calM}{{\mathcal M}}

\newcommand{\Ordo}{\mathrm{O}}

\newcommand{\bpi}{\boldsymbol\pi}

\newcommand{\im}{\operatorname{Im}}
\newcommand{\sign}{\operatorname{sgn}}
\newcommand{\Qop}{{\mathbf Q}}

\newcommand{\Jop}{{\mathbf J}}

\newcommand{\Hop}{{\mathbf H}}

\newcommand{\diff}{{\mathrm d}}
\newcommand{\imag}{{\mathrm i}}

\newcommand{\supp}{\operatorname{supp}}

\newcommand{\pv}{\operatorname{pv}}

\newcommand{\eps}{\varepsilon}

\begin{document}
%
\title[Fourier uniqueness sets and the Klein-Gordon equation]
{Fourier uniqueness sets and the Klein-Gordon equation}

\author[Hedenmalm]
{H\aa{}kan Hedenmalm}

\address{Hedenmalm: Department of Mathematics\\
The Royal Institute of Technology\\
S -- 100 44 Stockholm\\
SWEDEN}

\email{haakanh@math.kth.se}

\thanks{Research partially supported by the G\"oran Gustafsson Foundation 
(KVA) and by Vetenskapsr\aa{}det (VR)}









\begin{abstract} 
Building on ideas from \cite{HM}, we introduce (local) Fourier uniqueness 
sets for spaces of measures supported on a given curve in the plane. For 
the classical conic sections, the Fourier transform of the measure solves 
a second order partial diffeential equation. We focus mainly on the 
one-dimensional Klein-Gordon equation, which is associated with the hyperbola. 
We define the Hilbert transform for the hyperbola, and use it to introduce
a natural real Hardy space of absolutely continuous measures on the 
hyperbola. For that space of measures, we obtain several examples of (local) 
Fourier uniqueness sets. We also obtain examples of Fourier uniqueness sets 
in the context of all Borel measures on the curve. The proofs are based on 
the dynamics of Gauss-type maps combined with ideas from complex analysis. 
We also look at the Fourier uniqueness sets for one branch of the hyperbola, 
where the notion of defect becomes natural.  
\end{abstract}

\maketitle

\addtolength{\textheight}{2.2cm}







\section{Introduction}

\noindent\bf Heisenberg uniqueness pairs: variations on the theme. \rm
Let $\mu$ be a finite complex-valued Borel measure in the plane $\R^2$,
and associate to it the Fourier transform 
\[\widehat\mu(\xi):=\int_{\R^2}\e^{\pi\imag\langle x,\xi\rangle}\diff\mu(x),\]
where $x=(x_1,x_2)$ and $\xi=(\xi_1,\xi_2)$, with inner product
\[\langle x,\xi\rangle=x_1\xi_1+x_2\xi_2.\]
In \cite{HM}, the concept of a Heisenberg uniqueness pair (HUP) was introduced.
It is similar to the notion of (weakly) mutually annihilating pairs of Borel 
measurable sets having positive area measure, which appears, e.g., 
in the book by Havin and J\"oricke \cite{HJ}.
For $\Gamma\subset\R^2$ which is finite disjoint union of smooth curves in 
$\R^2$, let $\mathrm{M}(\Gamma)$ denote the Banach space of Banach space of 
complex-valued finite Borel measures in $\R^2$, supported on $\Gamma$. 
Moreover, let $\mathrm{AC}(\Gamma)$ denote the closed subspace of 
$\mathrm{M}(\Gamma)$ consisting of the measures that are absolutely 
continuous with respect to arc length measure on $\Gamma$.

\begin{defn}
Let $\Gamma$ be a finite disjoint union of smooth curves in $\R^2$, and let
$\mathrm{X}(\Gamma)$ be a linear subspace of $\mathrm{M}(\Gamma)$.
For a set $\Lambda\subset\R^2$, we say that {\em $\Lambda$ is a 
Fourier uniqueness set} for $\mathrm{X}(\Gamma)$ provided that
\[
\forall\mu\in\mathrm{X}(\Gamma):\quad
\widehat\mu|_{\Lambda}=0 \,\,\,\implies\,\,\,\mu=0.
\]
\end{defn}

Here, it is natural to require $\mathrm{X}(\Gamma)$ to be a norm closed 
subspace of $\mathrm{M}(\Gamma)$ (like, e.g., $\mathrm{AC}(\Gamma)$).
More generally, we could ask that $\mathrm{X}(\Gamma)$ is a {\em Banach
subspace} of $\mathrm{M}(\Gamma)$: this requires that 
$\mathrm{X}(\Gamma)$ is be equipped with a Banach space norm which makes
the injection mapping  $\mathrm{X}(\Gamma)\hookrightarrow\mathrm{M}(\Gamma)$  
continuous.

Following \cite{HM}, $(\Gamma,\Lambda)$ is a {\em Heisenberg uniqueness pair} 
if and only if $\Lambda$ is a Fourier uniqueness set for 
$\mathrm{AC}(\Gamma)$. 
The present concept offers the flexibility 
to consider more general spaces of measures.
We turn to the notion of a {\em defect}. 

\begin{defn}
Let $\Gamma$ be a finite disjoint union of smooth curves in $\R^2$, and let
$\mathrm{X}(\Gamma)$ be a linear subspace of $\mathrm{M}(\Gamma)$.
For a set $\Lambda\subset\R^2$, we say that $\Lambda$ is a Fourier
uniqueness set of {\em defect} $d$ for $\mathrm{X}(\Gamma)$ 
provided that the $\C$-linear space
\[
\big\{\mu\in\mathrm{X}(\Gamma):\,\,\widehat\mu|_{\Lambda}=0\big\}
\]
has dimension $d$.
\end{defn}

If we specialize to $\mathrm{X}(\Gamma)=\mathrm{AC}(\Gamma)$, the dual 
formulation is that $\Lambda$ is a Fourier uniqueness set of defect $d$ for 
$\mathrm{AC}(\Gamma)$ if and only if the weak-star closure  of the linear 
span of the complex exponentials
\[
e_\xi(x)=\e^{\pi\imag\langle x,\xi\rangle},\qquad \xi\in\Lambda,\,\,\,
x\in\Gamma,
\]
has codimension $d$ in $L^\infty(\Gamma)$. In line with the terminology of
\cite{HM}, we then say that $(\Gamma,\Lambda)$ is a {\em Heisenberg uniqueness
pair with defect $d$} (HUP${}_d$). The properties of the Fourier transform 
with respect to translation and multiplication by complex exponentials show 
that for all points $x^*,\xi^*\in\R^2$, we have
\[(\Gamma+\{x^*\},\Lambda+\{\xi^*\})\quad \text{is an HUP${}_d$}\quad
\Longleftrightarrow \quad (\Gamma,\Lambda)\quad \text{is an HUP${}_d$}.
\leqno{\text{(inv-1)}}\] 
Likewise, it is also straightforward to see that if $T:\R^2\to\R^2$ is an
invertible linear transformation with adjoint $T^*$, then
\[(T^{-1}(\Gamma),T^*(\Lambda))\quad \text{is an HUP}_{d}\quad
\Longleftrightarrow \quad (\Gamma,\Lambda)\quad \text{is an HUP}_d.
\leqno{\text{(inv-2)}}\] 
What is used here is a certain invariance of the space $\mathrm{AC}(\Gamma)$
under affine transformations of $\R^2$ as well as under multiplication by
complex exponentials; the analogous assertion would hold with more general 
classes of spaces $\mathrm{X}(\Gamma)$ provided those invariances remain valid.

We turn to the notion of local Fourier uniqueness sets. 

\begin{defn}
Let $\Gamma$ be a finite disjoint union of smooth curves in $\R^2$, and let
$\mathrm{X}(\Gamma)$ be a linear subspace of $\mathrm{M}(\Gamma)$. 
For two sets $\Lambda,\mathrm{K}\subset\R^2$, we say that 
$\Lambda$ is a {\em $\mathrm{K}$-local Fourier uniqueness set} 
(in short, FUS($\mathrm{K}$)) for $\mathrm{X}(\Gamma)$ provided that 
$\Lambda\subset\mathrm{K}$ and
\[
\forall\mu\in\mathrm{X}(\Gamma):\quad
\widehat\mu|_{\Lambda}=0 \,\,\,\implies\,\,\,\widehat\mu|_{\mathrm{K}}=0.
\]
\end{defn}

We generalize the notions of \cite{HM} by introducing local and strong 
Heisenberg uniqueness pairs.

\begin{defn}
If $\Lambda$ is a $\mathrm{K}$-local Fourier uniqueness set for 
$\mathrm{AC}(\Gamma)$, we say that $(\Gamma,\Lambda)$ is a 
{\em $\mathrm{K}$-local Heisenberg uniqueness pair} (in short, 
HUP$(\mathrm{K})$). If $\Lambda$ is a Fourier uniqueness set for 
$\mathrm{M}(\Gamma)$, we say that $(\Gamma,\Lambda)$ is a 
{\em strong Heisenberg uniqueness pair} (in short, 
SHUP). 
\label{defn-2}
\end{defn}

\begin{rem}
In terms of the familiar Zariski closure operation induced by the Fourier 
transforms of $\mathrm{X}(\Gamma)$, we are asking that the Zariski closure of 
$\Lambda$ should contain $\mathrm{K}$. We remark here that the space of 
Fourier transform from $\mathrm{X}(\Gamma)$ is generally speaking not an 
algebra, not even for $\mathrm{X}(\Gamma)=\mathrm{AC}(\Gamma)$; 
in case it is an algebra, the curve $\Gamma$ would necessarily be closed under 
addition. This means that we cannot expect to have a Zariski topology, 
although the closure operation is well-defined.
\end{rem}  

The dual formulation is that $(\Gamma,\Lambda)$ is an HUP($\mathrm{K}$)
if and only if the weak-star closure in $L^\infty(\Gamma)$ of the linear 
span of the functions $\{e_\xi(x)\}_{\xi\in\Lambda}$ contains all the 
functions $e_\xi(x)$ with $\xi\in\mathrm{K}$.

\begin{rem}
$(a)$ It is of course possible to mix the defect, locality, and strength 
notions, and talk about, e.g., HUP${}_d(\mathrm{K})$. We shall not need to 
do so in this presentation.  

\noindent $(b)$ If we write $\Lambda\prec_{\Gamma}\mathrm{K}$ to express that 
$(\Gamma,\Lambda)$ is a $K$-local Heisenberg uniqueness pair, we get a 
partial ordering. 
\end{rem}

We should understand the invariance properties of local Heisenberg uniqueness 
pairs with respect to affine transformations. As before, the properties of the 
Fourier transform with respect to translation and multiplication by complex 
exponentials show that for all points $x^*,\xi^*\in\R^2$, we have
\[(\Gamma+\{x^*\},\Lambda+\{\xi^*\})\quad \text{is an 
HUP}(\mathrm{K}+\{\xi^*\})\quad
\Longleftrightarrow \quad (\Gamma,\Lambda)\quad \text{is an HUP}(\mathrm{K}).
\leqno{\text{(inv-3)}}\] 
Likewise, it is also straightforward to see that if $T:\R^2\to\R^2$ is an
invertible linear transformation with adjoint $T^*$, then
\[(T^{-1}(\Gamma),T^*(\Lambda))\quad \text{is an HUP}(T^*(\mathrm{K}))\quad
\Longleftrightarrow \quad (\Gamma,\Lambda)\quad \text{is an HUP}(\mathrm{K}).
\leqno{\text{(inv-4)}}\] 
Again, what is used here is certain invariance of $\mathrm{AC}(\Gamma)$ under
affine transformations of $\R^2$ as well as under multiplication by complex 
exponentials; the analogous assertions would hold with more general classes
of spaces $\mathrm{X}(\Gamma)$ provided those invariance remain valid. 
\medskip

\noindent{\bf The Klein-Gordon equation.} 
In natural units, the Klein-Gordon equation reads
\begin{equation}
-\partial_t^2u+\Delta_xu=\EP u,
\label{eq-KG100}
\end{equation}
where $m>0$ is a constant (it is the mass of the particle), and 
\[
\Delta_x=\partial_{x_1}^2+\ldots+\partial_{x_d}^2
\]
is the $d$-dimensional Laplacian. We shall here restrict to the case of $d=1$,
one spatial dimension. So, our equation reads
\[
-\partial_t^2u+\partial_x^2u=\EP u.
\]
In terms of the (preferred) coordinates
\[
\xi_1:=x+t,\,\,\,\xi_2:=x-t,
\]
the Klein-Gordon equation reads
\[
\partial_{\xi_1}\partial_{\xi_2}u=\tfrac14\EP u.
\tag{KG} 
\]

\begin{rem}
Since $t^2-x^2=\xi_1\xi_2$, the space-like vectors (those vectors 
$(t,x)\in\R^2$ with $x^2-t^2>0$) correspond to the union of the first quadrant 
$\xi_1,\xi_2>0$ and the third quadrant $\xi_1,\xi_2<0$ in the 
$(\xi,\xi_2)$-plane). Likewise, the time-like vectors correspond to the
union of the second quadrant $\xi_1>0,\xi_2<0$ and the fourth quadrant
$\xi_1<0,\xi_2>0$.
\end{rem}
 
In the sequel, we will not need to talk about the time and space coordinates
$(t,x)$ as such. So, e.g., we are free to use the notation $x=(x_1,x_2)$ 
for the Fourier dual coordinate to $\xi=(\xi_1,\xi_2)$. 

Let $\calM(\R^2)$ denote the (Banach) space of all finite complex-valued
Borel measures in $\R^2$. We suppose that $u$ is the Fourier 
transform of a $\mu\in\calM(\R^{2})$:
\begin{equation}
u(\xi)=\hat\mu(\xi):=\int_{\R^{2}}
\e^{\pi\imag\langle x,\xi\rangle}
\diff\mu(x),\qquad\xi\in\R^2.
\label{eq-1.1'}
\end{equation}
The assumption that $u$ solves the Klein-Gordon equation (KG) means that
\[
\bigg(x_1x_2+\frac{\EP}{4\pi^2}\bigg)\diff\mu(x)=0
\]
as a measure on $\R^{2}$, which we see is the same as having
\begin{equation}
\supp\mu\subset\Gamma_m:=\bigg\{x\in\R^2:
\,x_1x_2=-\frac{\EP}{4\pi^2}\bigg\}.
\label{eq-1.2}
\end{equation}
The set $\Gamma_m$ is a hyperbola. We may use the $x_1$-axis to supply 
a global coordinate for $\Gamma_m$, and define a complex-valued finite Borel 
measure $\bpi_1\mu$ on $\R$ by putting
\begin{equation}
\bpi_1\mu(E)=\int_{E}\diff\bpi\mu(x_1):=\mu(E\times\R)=
\int_{E\times\R}\diff\mu(x).
\label{eq-bpi1}
\end{equation}
We shall at times refer to $\bpi_1\mu$ as the {\em compression} of $\mu$ to the
$x_1$-axis. It is easy to see that $\mu$ may be recovered from $\bpi_1\mu$; 
indeed,
\begin{equation}
u(\xi)=\hat\mu(\xi)=\int_{\R^\times}
\e^{\pi\imag[\xi_1t-m^2\xi_2/(4\pi^{2}t)]}
\diff\bpi_1\mu(t),\qquad \xi\in\R^2.
\label{eq-1.3}
\end{equation}
Here, we use the standard convention $\R^\times:=\R\setminus\{0\}$.
We note that $\mu$ is absolutely continuous with respect to arc length measure 
on $\Gamma_m$ if and only if $\bpi_1\mu$ is absolutely continuous with respect 
to (Lebesgue) length measure on $\R$. 
For positive reals $\alpha,\beta$, let $\Lambda_{\alpha,\beta}$
denote the lattice-cross 
\begin{equation}
\Lambda_{\alpha,\beta}:=(\alpha\Z\times\{0\})\cup(\{0\}\times\beta\Z),
\label{eq-LC}
\end{equation}
so that the spacing along the $\xi_1$-axis is $\alpha$, and along the 
$\xi_2$-axis it is $\beta$. In the recent paper \cite{HM}, Hedenmalm and 
Montes-Rodr\'{i}guez found the following.

\begin{thm} {\rm (Hedenmalm, Montes)} For positive reals $m,\alpha,\beta$, 
$(\Gamma_m,\Lambda_{\alpha,\beta})$ is a Heisenberg uniqueness pair if and
only if $\alpha\beta m^2\le 4\pi^2$.
\label{thm-1}
\end{thm}

Here, we consider possible generalizations. 

A first possibility is to ask what happens when 
$\alpha\beta m^2>4\pi^2$. By the above theorem, 
$(\Gamma_m,\Lambda_{\alpha,\beta})$ fails to be a Heisenberg uniqueness pair,
but it could still be a Heisenberg uniqueness pair with finite defect. 
This is, however, not the case, cf. \cite{CHM}.

Second, we may consider the union of the lattice-cross with a quadrant of the
plane. This way, we are able to produce Fourier uniqueness sets for 
$\mathrm{M}(\Gamma_m)$.  

Third, we may ask whether the portion of the lattice-cross which is in a 
given quadrant is a local Fourier uniqueness set. This problem is studied in 
Section \ref{sec-some}. 
Here, we need to consider a slightly smaller collection of 
measures than $\mathrm{AC}(\Gamma)$. 
If we combine the restrictions to quadrants with translations, we obtain 
interesting examples of Fourier uniqueness sets.  

Fourth, we may consider only one branch $\Gamma_m^+$ of the hyperbola 
$\Gamma_m$, and ask when the lattice-cross $\Lambda_{\alpha,\beta}$ is a 
Fourier uniqueness set for $\mathrm{AC}(\Gamma_m^+)$ in this setting. 
This problem can be understood in terms of when we have unique continuation
between the two branches of the hyperbola for measures in 
$\mathrm{AC}(\Gamma_m)$ whose Fourier transform vanishes on the lattice-cross.
\medskip

\noindent{\bf Acknowledgements.} I thank Alfonso Montes-Rodr\'\i{}guez
for several fruitful conversations. 

\section{ The Hilbert transform on the hyperbola}
\label{sec-Hilb}

\noindent{\bf Hilbert transforms.} 
The Hilbert transform $\Hop$ of a function $f\in L^1(\R)$ is
\[
\Hop[f](x):=\pv \frac1\pi \int_{\R}\frac{f(t)}{x-t}\diff t,\qquad x\in\R,
\]
wherever the integral makes sense. This may be thought of both as 
function in weak $L^1$, and as a distribution. We shall need to think of it
as a distribution. We easily extend the notion to measures:
for a finite complex-valued Borel measure $\nu$ on $\R$, we put
\[
\diff\Hop[\nu](x):=\bigg[\pv \frac1\pi\int_{\R}\frac{f(t)}{x-t}\diff t\bigg]
\diff x,
\]
where the notation suggests that we get a measure; this is just a formality,
as we generally expect only a distribution. As for the interpretation as
a weak $L^1$ function, we refer to the recent contribution \cite{PSZ} by
Poltoratski, Simon, and Zinchenko.

To simplify the notation, we restrict our attention to the hyperbola 
$\Gamma_m$ with $m=2\pi$, described by the equation $x_1x_2=-1$, and 
denote it by $\Gamma$ (dropping the subscript). Let us consider the 
following measure on $\Gamma$:
\[
\diff\lambda(x_1,x_2):=|x_1|^{-1}\diff\delta_{-1/x_1}(x_2)\diff x_1,
\]  
which has the symmetry property
\[
\diff\lambda(x_1,x_2)=\diff\lambda(x_2,x_1).
\]
For a finite Borel measure $\nu$ supported on $\Gamma$, we put
\[
\diff\Hop_\Gamma [\nu](x_1,x_2):=
\bigg[\pv\frac1\pi\int_\Gamma \frac{\diff\nu(y_1,y_2)}{x_1-y_1}\bigg]|x_1|\diff
\lambda(x_1,x_2),
\] 
which in general need not be a Borel measure, but rather can be interpreted as
a distribution supported on $\Gamma$. The way things are set up, $\bpi_1$
intertwines between $\Hop_\Gamma$ and $\Hop$:
\begin{equation}
\diff\bpi_1\Hop_\Gamma[\nu]=\diff\Hop[\bpi_1\nu].
\label{eq-Hilb2}
\end{equation}
After a moment's reflection we see that 
\[
\diff\Hop_\Gamma [\nu](x_1,x_2)=
\frac1\pi\sign(x_1)\nu(\R^2)\diff\lambda(x_1,x_2)
+
\bigg[\pv\frac1\pi\int_\Gamma \frac{\diff\nu(y_1,y_2)}{x_2-y_2}\bigg]|x_2|\diff
\lambda(x_1,x_2).
\]
Here, $\sign(t)$ is the sign of $t\in\R$ ($\sign(0)=0$, $\sign(t)=1$ for $t>0$,
and $\sign(t)=-1$ for $t<0$). If we let $\mathrm{AC}_0(\Gamma)$ denote the
codimension one subspace of $\mathrm{AC}(\Gamma)$ consisting of measures $\nu$
with $\nu(\R^2)=0$, we see that
\begin{equation*}
\diff\Hop_\Gamma [\nu](x_1,x_2)=
\bigg[\pv\frac1\pi\int_\Gamma \frac{\diff\nu(y_1,y_2)}{x_1-y_1}\bigg]|x_1|\diff
\lambda(x_1,x_2)=
\bigg[\pv\frac1\pi\int_\Gamma \frac{\diff\nu(y_1,y_2)}{x_2-y_2}\bigg]|x_2|\diff
\lambda(x_1,x_2),
\end{equation*}
which means that if $\bpi_2$ is the compression to the $x_2$-axis,
\[
\bpi_2\nu(E)=\int_E\diff\bpi_2\nu(x_2):=\int_{\R\times E}\diff\nu(x),
\]
then $\bpi_2$ intertwines $\Hop$ and $\Hop_\Gamma$ as well:
\begin{equation}
\diff\bpi_2\Hop_\Gamma[\nu]=\diff\Hop[\bpi_2\nu],\qquad \nu\in
\mathrm{AC}_0(\Gamma).
\label{eq-Hilb3}
\end{equation}
Since $\Hop_\Gamma$ relates to the Hilbert transform of the compression to
each of the two axes, it appears to be a rather natural operator.
We call it the {\em Hilbert transform on $\Gamma$}, and introduce the
{\em real $H^1$ space on $\Gamma$}, denoted $\mathrm{ACH}(\Gamma)$, which by
definition consists of those $\nu\in\mathrm{AC}_0(\Gamma)$ with 
$\Hop_\Gamma[\nu]\in\mathrm{AC}_0(\Gamma)$. Supplied with the
norm
\[
\|\nu\|_{\mathrm{ACH}(\Gamma)}:=\|\nu\|_{\mathrm{M}(\R^2)}+
\|\Hop_\Gamma[\nu]\|_{\mathrm{M}(\R^2)},\qquad \nu\in \mathrm{ACH}(\Gamma),
\]
it is a Banach space, and the injection $\mathrm{ACH}(\Gamma)\hookrightarrow
\mathrm{M}(\Gamma)$ is continuous, which makes $\mathrm{ACH}(\Gamma)$ a Banach
subspace of $\mathrm{M}(\Gamma)$ which is contained in $\mathrm{AC}_0(\Gamma)$.

\begin{rem}
In terms of the Fourier 
transform, we get from \eqref{eq-Hilb2} and \eqref{eq-Hilb3} that 
\begin{equation}
\forall \nu\in \mathrm{AC}_0(\Gamma):\,\,\,
\widehat{\Hop_\Gamma[\nu]}(\xi_1,0)=\imag\sign(\xi_1)\widehat\nu(\xi_1,0),\quad
\widehat{\Hop_\Gamma[\nu]}(0,\xi_2)=\imag\sign(\xi_2)\widehat\nu(0,\xi_2),
\label{eq-Hilb4}
\end{equation}
for all $\xi_1,\xi_2\in\R$. 
For $\nu\in\mathrm{AC}_0(\Gamma)$, the function $v:=\widehat \nu$ is 
continuous on $\R^2$, tends to $0$ at infinity (this is a consequence of the 
curvature of $\Gamma$), and has $v(0,0)=\widehat \nu(0,0)=\nu(\R^2)=0$. 
So it is immediate from \eqref{eq-Hilb4} that like $v=\widehat\nu$, the 
Fourier transform $v^*:=\widehat{\Hop_\Gamma[\nu]}$ solves the 
Klein-Gordon equation (KG) and has $v^*(\xi_1,\xi_2)=\imag\sign(\xi_1+\xi_2)
v(\xi_1,\xi_2)$ if $\xi_1\xi_2=0$, so $v^*(\xi_1,\xi_2)$ makes sense
as a continuous function on $\xi_1\xi_2=0$. Whether in 
general $v^*$ is automatically continuous throughout $\R^2$ is not so clear.
But if $\nu\in\mathrm{ACH}(\Gamma)$, there is of course no problem. 
\end{rem}

\section{Strong and weak Heisenberg uniqueness for the hyperbola}

\noindent{\bf Strong Heisenberg uniqueness for the hyperbola.} 
We recall the definition of strong Heisenberg uniqueness pairs (Definition 
\ref{defn-2}). 
First, we need some (standard) notation.
Let $\R_+$, $\R_-$ denote the sets of positive and negative reals, 
respectively, and put $\bar\R_{+}:=\R_+\cup\{0\}$, $\bar\R_{-}:=\R_-\cup\{0\}$.
We need the (standard) notion of a {\em Riesz set} $E_1\subset\R$: $E_1$ is
a Riesz set if $\mu\in\mathrm{M}(\R)$ and $\widehat\mu=0$ on $E_1$ implies that
$\mu$ is absolutely continuous. Here, we write
\[
\widehat\mu(\xi):=\int_\R\e^{\imag\pi\xi x}\diff\mu(x),\qquad \xi\in\R.
\]
By the well-known F. and M. Riesz theorem, any unbounded interval is a Riesz
set (see Proposition \ref{prop-0.5} below). This suggests the following 
definition.

\begin{defn}
Fix $m>0$. A set $E\subset\R^2$ is a {\em Riesz set for the hyperbola}
$\Gamma_m$ if every measure $\mu\in\mathrm{M}(\Gamma_m)$ with $\widehat\mu=0$
on $E$ is absolutely continuous with respect to arc length measure.
\end{defn}

\begin{rem}
We observe that given a Riesz set $E_1\subset\R$, the lifted sets
$E_1\times\{0\}$ and $\{0\}\times E_1$ are both Riesz sets for the hyperbola 
$\Gamma_m$.
\end{rem}

We consider the open quadrants 
\begin{equation}
\R^2_{++}:=\R_{+}\times\R_{+},\quad
\R^2_{--}:=\R_{-}\times\R_{-},\quad
\R^2_{+-}:=\R_{+}\times\R_{-},\quad
\R^2_{-+}:=\R_{-}\times\R_{+},
\label{eq-quad}
\end{equation}
and we write $\bar\R^2_{++},\bar\R^2_{--},\bar\R^2_{+-},\bar\R^2_{-+}$
for the corresponding closed quadrants.
The quadrants $\R^2_{++}$ and $\R^2_{--}$ are space-like, while $\R^2_{+-}$ and
$\R^2_{-+}$ are time-like.

\begin{thm}
Fix positive reals $\alpha,\beta,m$, and a Riesz set $E\subset\bar\R^2_{--}$
for the hyperbola. 
Then $\Lambda_{\alpha,\beta}^E:=\Lambda_{\alpha,\beta}\cup E$ 
is a Fourier uniqueness set for $\mathrm{M}(\Gamma_m)$ if and only if 
$\alpha\beta m^2\le 4\pi^2$. 
\label{thm-shup}
\end{thm}

In other words, for Riesz sets $E\subset\bar\R_-$, 
$(\Gamma_m,\Lambda_{\alpha,\beta}^{E})$ is a strong 
Heisenberg uniqueness pair provided that $\alpha\beta m^2\le 4\pi^2$.

\begin{rem}
The assertion remains the same if we replace the assumption 
$E\subset\bar\R^2_{--}$ by $E\subset\bar\R^2_{++}$.
\end{rem}
\medskip

\noindent{\bf Weak Heisenberg uniqueness for the hyperbola.} 
For general $m>0$, let $\mathrm{AC}_0(\Gamma_m)$ denote the subspace of
$\mathrm{AC}(\Gamma_m)$ consisting of measures $\nu$ with $\nu(\R^2)=0$.
A scaling argument allows us to define the Hilbert transform $\Hop_{\Gamma_m}$
on $\Gamma_m$ for general $m>0$, so that the analogue of \eqref{eq-Hilb4} 
holds:
\begin{equation}
\forall \nu\in \mathrm{AC}_0(\Gamma_m):\,\,\,
\widehat{\Hop_{\Gamma_m}[\nu]}(\xi_1,0)=
\imag\sign(\xi_1)\widehat\nu(\xi_1,0),\quad
\widehat{\Hop_{\Gamma_m}[\nu]}(0,\xi_2)=\imag\sign(\xi_2)\widehat\nu(0,\xi_2),
\label{eq-Hilb5}
\end{equation}
for all $\xi_1,\xi_2\in\R$. We define $\mathrm{ACH}(\Gamma_m)$ to be the
(dense) subspace of $\mathrm{AC}_0(\Gamma_m)$ of measures $\nu$ with 
$\Hop_{\Gamma_m}[\nu]\in\mathrm{AC}_0(\Gamma_m)$. This class of measures
is better-behaved, and it is quite natural to use it to define a slightly
bigger class of uniqueness sets.

\begin{defn}
The pair $(\Gamma_m,\Lambda)$ is a {\em weak Heisenberg uniqueness pair} if 
$\Lambda\subset\R^2$ is a Fourier uniqueness set for $\mathrm{ACH}(\Gamma_m)$.
Moreover, $(\Gamma_m,\Lambda)$ is a {\em $\mathrm{K}$-local weak Heisenberg
uniqueness pair} if $\Lambda$ (with $\Lambda\subset\mathrm{K}\subset\R^2$) 
is a $\mathrm{K}$-local Fourier uniqueness set for $\mathrm{ACH}(\Gamma_m)$.
\end{defn}

\begin{rem} $(a)$
As the terminology suggests, it is easier for $(\Gamma_m,\Lambda)$
to be a weak Heisenberg uniqueness pair than to be a Heisenberg uniqueness
pair. The same is true for the $\mathrm{K}$-local variant.

$(b)$ As $\mathrm{ACH}(\Gamma_m)\subset\mathrm{AC}_0(\Gamma_m)$ automatically,
we realize that for a set $\Lambda\subset\R^2$, we have the equivalence
(WHUP = weak Heisenberg uniqueness pair)
\[
(\Gamma_m,\Lambda)\,\,\,\text{is a WHUP} \,\,\,
\Longleftrightarrow\,\,\, 
(\Gamma_m,\Lambda\cup\{0\})\,\,\,\text{is a WHUP},
\] 
and, more generally, if $0\in\mathrm{K}$, we have
(WHUP(K) = $\mathrm{K}$-local weak Heisenberg uniqueness pair)
\[
(\Gamma_m,\Lambda)\,\,\,\text{is a WHUP(K)}\,\,\,\Longleftrightarrow\,\,\, 
(\Gamma_m,\Lambda\cup\{0\})\,\,\,\text{is a WHUP(K)}.
\] 
\end{rem}
Let $\mathrm{Q}$ be an open quadrant, i.e.,
\[
\mathrm{Q}\in\{\R^2_{++},\R^2_{--},\R^2_{+-},\R^2_{-+}\};
\]
we write $\bar{\mathrm{Q}}$ for the closure of $\mathrm{Q}$. 
The difference between space-like and time-like quarter-planes is made
obvious by the the following.

\begin{prop}
Fix a positive real $m$. Then:

\noindent{$(a)$}
For a time-like quarter-plane $\mathrm{Q}$, the boundary 
$\partial\mathrm{Q}$ is a Fourier uniqueness set for $\mathrm{M}(\Gamma_m)$.

\noindent{$(b)$} For a space-like quarter-plane $\mathrm{Q}$,  
the set $\bar{\mathrm{Q}}$ is not a Fourier uniqueness set for
$\mathrm{ACH}(\Gamma_m)$.
\label{prop-fund}
\end{prop}

We return to our lattice-cross $\Lambda_{\alpha,\beta}$ (see \eqref{eq-LC}),
and keep $\mathrm{Q}$ as an open quadrant.
Could it be that $\Lambda_{\alpha,\beta}\cap\bar{\mathrm{Q}}$ is a 
$\bar{\mathrm{Q}}$-local Fourier uniqueness set for $\mathrm{AC}(\Gamma_m)$,
or at least for $\mathrm{ACH}(\Gamma_m)$?
The answer to this question, as it turns out, depends on whether the quadrant
$\mathrm{Q}$ is space-like or time-like. For the time-like quarter-planes, 
there is no analogue of Theorem \ref{thm-1}, as can be seen from the following.

\begin{thm}
{\rm (Time-like quarter-planes)}
Fix $m,\alpha,\beta>0$. Then, for $\mathrm{Q}\in\{\R^2_{+-},\R^2_{-+}\}$,
the pair $\Lambda_{\alpha,\beta}\cap\bar{\mathrm{Q}}$ is not a 
$\bar{\mathrm{Q}}$-local Fourier uniqueness set for $\mathrm{ACH}(\Gamma_m)$.
In particular, $\Lambda_{\alpha,\beta}\cap\bar{\mathrm{Q}}$ is not 
a Fourier uniqueness pair for $\mathrm{ACH}(\Gamma_m)$.
%
%
\label{thm-2}
\end{thm} 

As regards the space-like quarter-planes, there is indeed an analogue.

\begin{thm}
{\rm (Space-like quarter-planes)}
Fix positive reals $m,\alpha,\beta$.
Then, for $\mathrm{Q}\in\{\R^2_{++},\R^2_{--}\}$, 
$\Lambda_{\alpha,\beta}\cap\bar{\mathrm{Q}}$ is a $\bar{\mathrm{Q}}$-local 
Fourier uniqueness set for $\mathrm{ACH}(\Gamma_m)$ if and only if 
$\alpha\beta m^2\le 4\pi^2$.
\label{thm-3}
\end{thm}


\begin{rem}
It is not known whether for $\alpha\beta m^2\le 4\pi^2$ the set
$\Lambda_{\alpha,\beta}\cap\bar{\mathrm{Q}}$ is a $\bar{\mathrm{Q}}$-local 
Fourier uniqueness set for $\mathrm{AC}(\Gamma_m)$. 
This problem appears rather challenging, and it has an attractive 
reformulation (see Problems \ref{prob-1} and \ref{prob-2}).
\label{rem-1}
\end{rem}
\medskip

\noindent{\bf A family of weak Heisenberg uniqueness pairs for the hyperbola.}
For a point $\xi^0=(\xi^0_1,\xi^0_2)\in\R^2$, we consider the distorted 
lattice-cross
\[
\Lambda_{\alpha,\beta}^{\langle\xi^0\rangle}:=(\Lambda_{\alpha,\beta}\cap
\bar\R^2_{--})\cup((\Lambda_{\alpha,\beta}\cap
\bar\R^2_{++})+\{\xi^0\}),
\]
which has the general appearance of a ``slanted lattice-waist'' if 
$\xi^0_1>0$ or $\xi^0_2>0$. 

\begin{thm}
Suppose $\alpha,\beta,m$ are positive reals, with $\alpha\beta m^2
\le 4\pi^2$. Then $(\Gamma_m,\Lambda_{\alpha,\beta}^{\langle\xi^0\rangle})$ 
is a weak Heisenberg uniqueness pair if and only if 
$\min\{\xi^0_1,\xi^0_2\}<0$ or $\xi^0=(0,0)$.
\label{thm-3.8}
\end{thm}

\begin{rem}
For $\alpha\beta m^2>4\pi^2$, we do not know what happens, but we suspect
that $\Lambda_{\alpha,\beta}^{\langle\xi^0\rangle}$ is not a weak Heisenberg 
uniqueness pair, independently of the location of $\xi^0$. 
\end{rem}

\section{Heisenberg uniqueness for one branch of the hyperbola}

\noindent{\bf The branches of the hyperbola.}
The hyperbola $\Gamma_m$ naturally splits into two connectivity components:
\begin{equation}
\Gamma_m^+:=\bigg\{x\in\R^2:
\,x_1x_2=-\frac{\EP}{4\pi^2},\,\,x_1>0\bigg\},\quad
\Gamma_m^-:=\bigg\{x\in\R^2:
\,x_1x_2=-\frac{\EP}{4\pi^2},\,\,x_1<0\bigg\}.
\label{eq-1.5}
\end{equation}
In light of Theorem \ref{thm-1}, it is natural to ask what happens if we
replace $\Gamma_m$ by one of $\Gamma_m^+,\Gamma_m^-$. In view of the 
invariance property (inv-2), it suffices to treat $\Gamma_m^+$.

\begin{thm}
For positive reals $m,\alpha,\beta$, $(\Gamma_m^+,\Lambda_{\alpha,\beta})$ is a
Heisenberg uniqueness pair if and only if $\alpha\beta m^2<16\pi^2$. Moreover,
for $\alpha\beta m^2=16\pi^2$, $(\Gamma_m^+,\Lambda_{\alpha,\beta})$ is a
Heisenberg uniqueness pair with defect $1$. 
\label{thm-onebranch}
\end{thm}

\begin{rem} 
We suspect that for $\alpha\beta m^2>16\pi^2$, 
$(\Gamma_m^+,\Lambda_{\alpha,\beta})$ is not a Heisenberg uniqueness pair 
with a finite defect $d$ (i.e., the defect should be infinite). Cf. Remark
\ref{rem-g1}.
\end{rem}

In the critical case $\alpha\beta m^2=16\pi^2$, we can get rid of the defect
by adding a point on the cross which does not lie on the lattice-cross.

\begin{cor}
Suppose $m,\alpha,\beta$ are positive reals with $\alpha\beta m^2=16\pi^2$.
Pick a point $\xi^0\in(\R\times\{0\})\times(\{0\}\times\R)$ on the cross,
and put $\Lambda^0_{\alpha,\beta}:=\Lambda_{\alpha,\beta}\cup\{\xi^0\}$.
Then $(\Gamma_m^+,\Lambda^0_{\alpha,\beta})$ is a Heisenberg uniqueness 
pair if and only if $\xi^0\not\in\Lambda_{\alpha,\beta}$.
\label{cor-onebranch}
\end{cor}

\section{Elements of Hardy space theory}
\label{sec-hardy}

\noindent{\bf Hardy spaces.} We shall need certain subspaces of $L^1(\R)$
and $L^\infty(\R)$. If $f$ is in $L^1(\R)$ or in $L^\infty(\R)$, we define
its Poisson extension to the upper half-plane
\[
\C_+:=\{z\in\C:\,\im z>0\}
\]
by the formula
\[
\poiss_{\C_+}[f](z):=\frac{\im z}{\pi}\int_\R\frac{f(t)}{|z-t|^2}\,\diff t,
\qquad z\in\C_+.
\]
The function $\poiss_{\C_+}[f]$ is harmonic in $\C_+$, and its boundary values 
are those of $f$ in the natural sense. It is standard to identify the function 
$f$ with its Poisson extension. We say that $f\in H^1_+(\R)$ if $f\in L^1(\R)$ 
and $\poiss_{\C_+}[f]$ is holomorphic. Likewise, we say that 
$f\in H^\infty_+(\R)$ if $f\in L^\infty(\R)$ and $\poiss_{\C_+}[f]$ is 
holomorphic. Analogously, if $f\in L^1(\R)$ and $\poiss_{\C_+}[f]$ is 
conjugate holomorphic (this means that the complex conjugate is holomorphic), 
we say that $f\in H^\infty_-(\R)$, while if $f\in L^\infty(\R)$ and 
$\poiss_{\C_+}[f]$ is conjugate holomorphic, we write 
$f\in H^\infty_-(\R)$. Clearly, $f\in H^1_-(\R)$ if and only if its
complex conjugate is in $H^1_+(\R)$, and the same goes for $H^\infty_-(\R)$
and $H^\infty_+(\R)$.

We shall use the following bilinear form on $\R$:
\[
\langle f,F\rangle_\R:=\int_\R f(t)F(t)\diff t,
\]
whenever it is well-defined. 
We shall frequently need the following well-known characterization of 
$H^1_+(\R)$. 

\begin{prop}
\label{prop-0}
Let us agree to write $e_\tau(t):=\e^{\imag\pi\tau t}$. Then the following
are equivalent for a function $f\in L^1(\R)$: $(a)$ $f\in H^1_+(\R)$, and
$(b)$  $\langle f,e_\tau\rangle_\R=0$ for all $\tau>0$.
\end{prop}

The following result is also standard.

\begin{prop}
\label{prop-0'}
\noindent{$(a)$} If $f\in H^1_+(\R)$ and $F\in H^\infty_+(\R)$, then 
$Ff\in H^1_+(\R)$, and $\langle f,F\rangle_\R=0$. 

\noindent{$(b)$} If $f\in L^1(\R)$, then $f\in H^1_+(\R)$ if and only if
$\langle f,F\rangle_\R=0$ for all $F\in H^\infty_+(\R)$. 
\end{prop}

We need also the next result, attributed to F. and M. Riesz.

\begin{prop}
\label{prop-0.5}
Suppose $\mu$ is a complex-valued finite Borel measure on $\R$. If
\[
\forall \tau>0:\quad \int_\R \e^{\imag\pi\tau t}\diff\mu(t)=0,
\]
then $\mu$ is absolutely continuous, and $\diff\mu(t)=f(t)\diff t$, where
$f\in H^1_+(\R)$.
\end{prop}
\medskip

\noindent{\bf Applications of Hardy space methods.}
We now show (see Proposition \ref{prop-2} below) that the linear span of the 
functions 
\[
t\mapsto\e^{\pi\imag\xi_1t},\quad t\mapsto \e^{\pi\imag\epsilon\xi_2/t},
\qquad\xi=(\xi_1,\xi_2)\in\bar\R^2_{++},
\]
is weak-star dense in $L^\infty(\R)$. 

\begin{prop}
Let $\nu$ be a complex-valued finite Borel measure on $\R$. If
\[
\int_\R\e^{\pi\imag\xi_1 t}\diff\nu(t)=
\int_\R\e^{\pi\imag\xi_2/t}\diff\nu(t)=0,\quad\text{for all}\,\,\,
\xi_1,\xi_2>0,
\]
then $\nu=0$.
\label{prop-2}
\end{prop}

\begin{proof}
By Proposition \ref{prop-0.5}, the assumptions entail
that $\diff\nu(t)=f_1(t)\diff t$ and $\diff\nu(-1/t)=f_2(-t)\diff t$, where
$f_j\in H^1_+(\R)$, $j=1,2$. 
By equating the two ways to represent $\diff\nu$, we see that
\[
f_1(t)=t^{-2}f_2(1/t),\qquad t\in\R,
\]
so that $f_1$ has an analytic pseudocontinuation to the lower half-plane
(for $\im t>0$, $\im t^{-1}<0$). The pseudocontinuation is of course a genuine
holomorphic continuation to $\C^\times:=\C\setminus\{0\}$ (we can use, e.g., 
Morera's theorem). In terms of 
$g_j(t):=tf_j(t)$, for $j=1,2$, the above relation reads $g_1(t)=g_2(1/t)$.
The functions $g_j$, $j=1,2$, extend holomorphically to $\C^\times$, and have
the estimate
\[
|g_j(t)|\le\|\nu\|\frac{|t|}{|\im t|},\qquad t\in\C\setminus\R,\,\,\,j=1,2,
\]
Using the theory around the $\log\log$ theorem (attributed to Levinson,
Sj\"oberg, Carleman, Beurling; see e.g. \cite{PK}, pp. 374--383, also 
\cite{BH}) it is not difficult to show that such functions
$g_j$, $j=1,2$, must be constant. But then the constant must be $0$, for
otherwise, $f_j$, $j=1,2$, would not be in $H^1_+(\R)$.
\end{proof}

\begin{proof}[Proof of Proposition \ref{prop-fund}]
We first consider time-like quarter-planes $\mathrm{Q}\in\{\R^2_{+-},
\R^2_{-+}\}$. In both cases, the problem boils down to Proposition 
\ref{prop-2}, which settles the issue.

We turn to space-like quarter-planes $\mathrm{Q}\in\{\R^2_{++},\R^2_{--}\}$.
Here, the matter is settled by Proposition \ref{prop-0'}.
The non-trivial measures $\mu\in\mathrm{ACH}(\Gamma_m)$ whose Fourier transform
vanishes on $\bar{\mathrm{Q}}$ have compressions to the $x_1$-axis of the form
$f(t)\diff t$, where $f\in H^1_+(\R)$ or $f\in H^1_-(\R)$ (which of the two it
is depends on whether $\mathrm{Q}$ is $\R^2_{++}$ or $\R^2_{--}$).  
\end{proof}

We next show (see Proposition \ref{prop-3} below) that the linear span of 
the functions 
\[
t\mapsto\e^{\pi\imag\xi_1t},\quad t\mapsto \e^{\pi\imag\xi_2/t},
\qquad\xi=(\xi_1,\xi_2)\in\bar\R^2_{+-},
\]
is weak-star dense in $H^\infty_+(\R)$. 

\begin{prop}
Let $\nu$ be a complex finite absolutely continuous measure on $\R$. Then 
\[
\int_\R\e^{\pi\imag\xi_1 t}\diff\nu(t)=
\int_\R\e^{-\pi\imag\xi_2/t}\diff\nu(t)=0,\quad\text{for all}\,\,\,
\xi_1,\xi_2>0,
\]
if and only if $\diff\nu(t)=f(t)\diff t$ where $f\in H^1_+(\R)$.
\label{prop-3}
\end{prop}

\begin{proof}
The assertion is immediate from Propositions \ref{prop-0} and \ref{prop-0'},
once it is observed that the function $F(t)=\e^{\pi\imag\xi_2/t}$ is in 
$H^\infty(\R)$ for $\xi_2<0$.  
\end{proof}


\medskip

\noindent{\bf The predual of $H^\infty$ on the line.} It is well-known and
can be seen from Proposition \ref{prop-0'} that the predual of 
$H^\infty_+(\R)$ is the quotient space $L^1(\R)/H^1_+(\R)$ with respect to 
the standard bilinear form $\langle\cdot,\cdot\rangle_\R$.
For $f\in H^1_+(\R)$ and $F\in H^\infty_+(\R)$, we have -- by Proposition
\ref{prop-0'} -- $\langle f,F\rangle_\R=0$, which is why we need to mod out 
with respect to $H^1_+(\R)$ in the predual. 
\medskip

\noindent{\bf The dual of $H^1$ on the line.} If we put
\[
H^1_{\mathrm{real}}(\R):=H^1_+(\R)\oplus H^1_-(\R),
\]
and supply this space with the natural norm; as $H^1_+(\R)\cap H^1_-(\R)
=\{0\}$, this is just the sum of the two norms:
\[
\|f_1+f_2\|_{H^1_{\mathrm{real}}(\R)}:=\|f_1\|_{H^1_+(\R)}+\|f_2\|_{H^1_-(\R)},
\qquad 
f_1\in H^1_+(\R),\,\, f_2\in H^1_-(\R).
\]
The {\em Cauchy projection} 
\[
\proj_+: H^1_{\mathrm{real}}(\R)\to H^1_+(\R),\quad 
\proj_+[f_1+f_2]:=f_1 \,\,\,\text{for}\,\,f_1\in H^1_+(\R),\,\, 
f_2\in H^1_-(\R),
\]  
is a thus norm contraction. It is related to the Hilbert transform $\Hop$:
\[
\proj_+f=\tfrac12(f+\imag\Hop[f]),\qquad f\in H^1_{\mathrm{real}}(\R).
\] 
The space $H^1_{\mathrm{real}}(\R)$ is a Banach space, and a dense (Banach) 
subspace of
\[
L^1_0(\R):=\{f\in L^1(\R):\,\,\langle f,1\rangle_\R=0\}.
\]
Actually, the space $H^1_{\mathrm{real}}(\R)$ has an alternative 
characterization in terms of the Hilbert transform:
\[
H^1_{\mathrm{real}}(\R)=\big\{f\in L^1_0(\R):\,\,\Hop[f]\in L^1_0(\R)\big\}.
\]
The dual space of $L^1_0(\R)$ is $L^\infty(\R)/\{\text{constants}\}$. 
The dual space of $H^1_{\mathrm{real}}(\R)$ is $\mathrm{BMO}(\R)$, which is 
understood as the space of functions with {\em bounded mean oscillation}, 
modulo the constants. The Cauchy projection also acts on the dual side:
\[
\proj_+:\text{BMO}(\R)\to \text{BMOA}_+(\R),
\]
where $\text{BMOA}_+(\R)$ is the subspace of $\text{BMO}(\R)$ which is dual 
to $H^1_-(\R)$ with respect to $\langle\cdot,\cdot\rangle_\R$. 

\medskip

\noindent{\bf The predual of $H^\infty$ on the unit circle.} 
We also need Hardy spaces in the context of the unit circle (or the unit 
disk, if we talk about the harmonic extension). A function in $L^1(\T)$ 
($\T$ is the unit circle) has norm 
\[
\|f\|_{L^1(\T)}:=\int_{-\pi}^{\pi}|f(\e^{\imag t})|\,\frac{\diff t}{2\pi},
\]
and we use the standard bilinear form
\[
\langle f,g\rangle_\T:=\int_{-\pi}^{\pi} f(\e^{\imag t})g(\e^{\imag t})\,
\frac{\diff t}{2\pi},\qquad f\in L^1(\T),\,\,\,
g\in L^\infty(\T).
\]
The Poisson extension to the unit disk $\D$ of $f\in L^1(\T)$ is given by the 
formula
\[
\poiss_\D f(z):=\int_{-\pi}^{\pi}\frac{1-|z|^2}{|1-z\e^{-\imag t}|^2}\,
f(\e^{\imag t})\,\frac{\diff t}{2\pi},\qquad z\in\D.
\]
If $f\in L^1(\T)$ and $\poiss_\D f$ is holomorphic in $\D$, we write 
$f\in H^1_+(\T)$. If, in addition, $\poiss_\D f(0)=0$, we write 
$f\in H^1_{+,0}(\T)$.  
We frequently identify functions on the unit circle $\T$ with their harmonic
extensions to $\D$. If $H^1(\D),H^1_0(\D)$ are defined as the spaces of
such extensions of boundary functions, we thus identify $H^1_+(\T)
\cong H^1(\D)$, $H^1_{+,0}(\T)\cong H^1_0(\D)$. In a similar fashion, 
$H^\infty_+(\T)\cong H^\infty(\D)$.
It is well-known that with respect to the standard bilinear form, the predual 
of $H^\infty_+(\T)$ may be identified with $L^1(\T)/H^1_{+,0}(\T)$. 
\medskip

%

\noindent{\bf Periodic Hardy spaces and the exponential mapping.}
Let $L^\infty(\R/2\Z)$ consist of those $f\in L^\infty(\R)$ which are
$2$-periodic: $f(x+2)=f(x)$. Similarly, we let $H^\infty_+(\R/2\Z)$ consist 
the $2$-periodic functions in $H^\infty_+(\R)$. 
The exponential mapping $x\mapsto\e^{\imag\pi x}$ provides an identification 
$\R/2\Z\cong\T$, and the upper half space $\C_+$ modulo $2\Z$ corresponds to 
the punctured disk $\D\setminus\{0\}$. The results for
the unit circle therefore carry over in a natural fashion to the $2$-periodic
setting. We let $L^1(\R/2\Z)$ denote the space of locally integrable 
$2$-periodic functions on $\R$, supplied with the Banach space norm
\[
\|f\|_{L^1(\R/2\Z)}:=\int_{[-1,1]}|f(x)|\diff x.
\]
We let $H^1_+(\R/2\Z)$ denote the subspace of $L^1_(\R/2\Z)$ consisting
of functions whose Poisson extension to the upper half plane $\C_+$ are
holomorphic. The holomorphic extension is then automatically $2$-periodic, 
and if, for $f\in H^1_+(\R/2\Z)$, the holomorphic extension (also denoted by 
$f$) has $f(z)\to0$ as $\im z\to+\infty$, we write $f\in H^1_{+,0}(\R/2\Z)$.
Via the exponential mapping, $H^1_{+}(\R/2\Z)$ corresponds to $H^1(\T)$, and
$H^1_{+,0}(\R/2\Z)$ to $H^1_{+,0}(\T)$. By carrying over the results available 
in the setting of the circle $\T$, we see that with respect to the bilinear 
form
\[
\langle f,g\rangle_{[-1,1]}:=\int_{[-1,1]}f(x)g(x)\diff x,\qquad f\in 
L^1_{(2)}(\R),\,\,\, g\in L^\infty_{(2)}(\R),
\]  
$H^1_{+,0}(\R/2\Z)$ is the pre-annilator of $H^\infty_{+}(\R/2\Z)$, and 
we may identify
\[
\big[L^1(\R/2\Z)/H^1_{+,0}(\R/2\Z)\big]^*=H^\infty_{+}(\R/2\Z).
\]

\section{Some reformulations and proofs}
\label{sec-some}

\noindent{\bf Strong Heisenberg uniqueness for the hyperbola.}
We may now supply the proof of  Theorem \ref{thm-shup}.

\begin{proof}[Proof of Theorem \ref{thm-shup}.]
We consider $\mu\in\mathrm{M}(\Gamma_m)$ with $\widehat\mu=0$ on 
$\Lambda_{\alpha,\beta}^{E}$. The assumption that $E$ is a Riesz set
for $\Gamma_m$ entails that the $\mu$ is absolutely continuous with respect 
to arc length measure on $\Gamma_m$. 
The main theorem of \cite{HM} -- based on the dynamics of the 
Gauss-type map $t\mapsto -\beta/t$ modulo $2$ on the interval $]-1,1]$ -- 
shows that for $\alpha\beta m^2\le4\pi^2$, the assumption that $\widehat\mu=0$ 
on $\Lambda_{\alpha,\beta}$ implies that $\mu=0$ identically. 
If we use that $E\subset\bar\R_-$, we can adapt the counterexample from 
\cite{HM} -- involving harmonic extensions -- to construct non-trivial measures
$\mu\in\mathrm{ACH}(\Gamma_m)$ with $\widehat\mu=0$ on 
$\Lambda_{\alpha,\beta}^E$ in case $\alpha\beta m^2>4\pi^2$. 
The proof is complete.
\end{proof}
\medskip

\noindent{\bf The dual formulation for time-like quarter-planes.}
Theorem \ref{thm-2} deals with the time-like quarter-planes 
$\mathrm{Q}\in\{\R^2_{+-},\R^2_{-+}\}$. 
If we take the invariance (inv-4) into account,
with $T$ as the reflection in the origin $(x_1,x_2)\mapsto(-x_1,-x_2)$, we 
realize that it suffices to consider $\mathrm{Q}=\R^2_{++}$. The dual 
formulation of the theorem runs as follows. 
{\em For all triples $\alpha,\beta,m>0$, the linear span of the functions
\[
\e^{\pi\imag\alpha jt},\,\,\,\e^{\imag\beta m^2k/(4\pi t)},\qquad
j,k=0,1,2,\ldots,
\]
fails to be weak-star dense in $L^\infty(\R)$.} 
By a scaling argument, we may assume that
\[
\alpha=1,\quad m=2\pi,
\]
so that we are dealing with the linear span of 
\[
\e^{\pi\imag jt},\,\,\,\e^{\pi\beta\imag k/t},\qquad
j,k=0,1,2,\ldots.
\]
This dual formulation of course requires that have Proposition \ref{prop-fund} 
at our disposal. Actually, Proposition \ref{prop-fund} may be deduced in a 
straightforward fashion from Propositions \ref{prop-2} and \ref{prop-3}. 
We leave the necessary details to the reader. This allows us to proceed with 
the proof of Theorem \ref{thm-2}. 

\begin{proof}[Proof of Theorem \ref{thm-2}]
First, we note that the functions $\e^{\pi\imag jt}$ belong to $H^\infty_+(\R)$
for $j=0,1,2,\ldots$, while the functions $\e^{\pi\beta\imag k/t}$ instead 
belong to $H^\infty_-(\R)$ for $k=0,1,2,\ldots$. This means that
the spanning vectors live in rather different subspaces and have no chance 
to span $\mathrm{BMO}(\R)$ even after weak-star closure. To make this more 
concrete, we pick a point $z_0\in\C_+$ in the upper half-plane and consider the
function 
\[
f_{z_0}(t):=\frac{1}{t-z_0}-\frac{1}{t-2-z_0},\qquad t\in\R.
\]
Clearly, $f_{z_0}(t)=\Ordo(t^{-2})$ as $|t|\to+\infty$, and so 
$f_{z_0}\in L^1(\R)$. Actually, we have $f_{z_0}\in H^1_-(\R)
\subset H^1_{\mathrm{real}}(\R)$.
We may use the calculus of residue to obtain that
\[
\int_\R f_{z_0}(t)\,\e^{\pi\imag jt}\diff t=2\pi\imag(\e^{\pi\imag jz_0}
-\e^{\pi\imag j(z_0+2)})=0,\qquad j=0,1,2,\ldots.
\]
Next, we may show that
\[
\int_\R f_{z_0}(t)\,\e^{\pi\beta\imag k/t}\diff t=0,\qquad k=0,1,2,\ldots,
\]
by appealing to Proposition \ref{prop-0'} (we will need to take complex 
conjugates if we work in the setting of the upper half-plane). 
So, for each $z_0\in\C_+$, $f_{z_0}$ annihilates the subspace, which 
consequently cannot be weak-star dense.
\end{proof}

\begin{rem}
The argument of the proof of theorem \ref{thm-2} actually shows that
the weak-star closure of the subspace spanned by $\e^{\pi\imag jt},
\e^{\pi\beta\imag k/t}$, for $j,k=0,1,2,\ldots$, has infinite codimension
in $\mathrm{BMO}(\R)$.
\end{rem}

\noindent{\bf The dual formulation for space-like quarter-planes.}
Theorem \ref{thm-3} and the open problem mentioned in Remark \ref{rem-1} deal 
with the space-like quarter-planes $\mathrm{Q}
\in\{\R^2_{++},\R^2_{--}\}$. If we take the invariance (inv-4) into account,
with $T$ as the inversion $(x_1,x_2)\mapsto(-x_1,x_2)$, we realize that it
suffices to consider $\mathrm{Q}=\R^2_{+-}$. The dual formulation of the 
theorem runs as follows. {\em For all triples of positive numbers 
$\alpha,\beta,m$, the linear span of the functions
\[
\e^{\pi\imag\alpha jt},\,\,\,\e^{-\imag\beta m^2k/(4\pi t)},\qquad
j,k=0,1,2,\ldots,
\]
(taken modulo the constants) is weak-star weak-star dense in 
$\mathrm {BMOA}_+(\R)$ if and only if $\alpha\beta m^2\le4\pi^2$.}
Alternatively, given $\mu\in\mathrm{ACH}(\R)$, we consider its compression
to the $x_1$-axis $\bpi_1\mu$, which has $\diff\bpi_1\mu(t)=f(t)\diff t$,
where $f\in H^1_{\mathrm{real}}(\R)$. We need to show that
\[
\big\langle t\mapsto\e^{\pi\imag\alpha jt},f\big\rangle_\R
=\big\langle t\mapsto \e^{-\imag\beta m^2k/(4\pi t)},f\rangle_\R=0,\qquad
j,k=0,1,2,\ldots,
\]
entails that $f\in H^1_+(\R)$ if and only if $\alpha\beta m^2\le4\pi^2$.

\begin{proof}[Proof of Theorem \ref{thm-3}]
The necessity of the condition $\alpha\beta m^2\le4\pi^2$ is just as in 
\cite{HM}, so we focus on the sufficiency.
We split $f=f_1+f_2$, where $f_1\in H^1_+(\R)$ and $f_2\in H^1_-(\R)$. Now, if
we apply Proposition \ref{prop-3} to $f_1$, we conclude that
\begin{equation}
\big\langle t\mapsto\e^{\pi\imag\alpha jt},f_2\big\rangle_\R
=\big\langle t\mapsto \e^{-\imag\beta m^2k/(4\pi t)},f_2\rangle_\R=0,\qquad
j,k=0,1,2,\ldots.
\label{eq-f2}
\end{equation}
Next, as $f_2\in H^1_-(\R)$, \eqref{eq-f2} actually holds for all $j,k\in\Z$. 
This puts us in the setting of \cite{HM}, and we find that $f_2=0$. The 
claim $f\in H^1_+(\R)$ follows.
\end{proof}
\medskip

\noindent{\bf An open problem for space-like quarter-planes.}  
We now turn to the open problem mentioned in Remark \ref{rem-1}.
By a scaling argument, we may assume that
\[
\alpha=1,\quad m=2\pi,
\]
so that we are dealing with the linear span of 
\[
\e^{\imag\pi jt},\,\,\,\e^{-\imag\pi\beta k/t},\qquad
j,k=0,1,2,\ldots.
\]
The issue at hand is whether this linear span is weak-star dense in 
$H^\infty(\R)$ for $\beta\le1$.  So, if $f\in L^1(\R)$ has 
\[
\big\langle t\mapsto\e^{\imag\pi jt},f\big\rangle_\R
=\big\langle t\mapsto\e^{-\imag\pi\beta k/t},f\big\rangle_\R=0,\qquad
j,k=0,1,2,\ldots
\]
may we then conclude (for $\beta\le1$) that $f\in H^1_+(\R)$?
For $j=0,1,2,\ldots$, the functions $t\mapsto\e^{\imag\pi jt}$
belong to $H^\infty_+(\R)$, and they are $2$-periodic: $\e^{\imag\pi j(t+2)}
=\e^{\imag\pi jt}$. From the well-known theory of Fourier series we obtain 
that the linear span of these functions $t\mapsto\e^{\imag\pi jt}$, where
$j=0,1,2,\ldots$, is weak-star dense in $H^\infty_{+}(\R/2\Z)$, 
the subspace of
$2$-periodic $H^\infty(\R)$ functions. As for the remaining spanning vectors
$\e^{-\imag\pi\beta k/t}$ a similar argument shows that their linear span
is weak-star dense in $H^\infty_+(\R/\langle\beta\rangle)$. Here, $g\in
H^\infty_+(\R/\langle\beta\rangle)$ if and only if $g\in H^\infty(\R)$ 
has $\{t\mapsto g(-\beta/t)\}\in H^\infty(\R/2\Z)$. In other words, 
$g\in H^\infty_+(\R/\langle\beta\rangle)$ means that $g\in H^\infty_+(\R)$ 
has the ``M\"obius periodicity''
\[
g\bigg(\frac{\beta t}{\beta-2t}\bigg)=g(t).
\]
We reformulate the problem in terms of these subspaces of $H^\infty_+(\R)$.

\begin{prob}
Is the sum $H^\infty_+(\R/2\Z)+H^\infty_+(\R/\langle\beta\rangle)$  
weak-star dense in $H^\infty(\R)$ for $\beta\le1$?
\label{prob-1}
\end{prob}

We reformulate this problem in terms of the periodization  operator
$\Qop_2: L^1(\R)\to L^1(\R/2\Z)$:
\begin{equation}
\Qop_2 f(x):=\sum_{j\in\Z}f(x+2j).
\label{eq-sum1}
\end{equation}
We first look at what it means for a function $f\in L^1(\R)$ that
\begin{equation}
\langle f,g\rangle_\R=0\quad \text{for all}\,\,\, g\in H^\infty_+(\R/2\Z).
\label{eq-annih1}
\end{equation}
For $f\in L^1(\R)$ and $g\in L^\infty(\R/2\Z)$, we see that
\begin{equation*}
\langle f,g\rangle_\R=\int_\R f(x)g(x)\diff x=\sum_{j\in\Z}\int_{[2j-1,2j+1]}
f(x)g(x)\diff x=\int_{[-1,1]}\Qop_2 f(x)g(x)\diff x
=\langle\Qop_2 f,g\rangle_{[-1,1]},
\end{equation*}
Via the exponential map $z\mapsto \e^{\imag\pi z}$ the space $H^\infty_+(\R
/2\Z)$ can be identified with $H^\infty_+(\T)$, and in view of the 
identification of the pre-annihilator of $H^\infty_+(\R/2\Z)$,
we find that \eqref{eq-annih1} is equivalent to having 
\begin{equation}
\Qop_2 f\in H^1_{+,0}(\R/2\Z).
\label{eq-annih2}
\end{equation}
We turn to the interpretation of
\begin{equation}
\langle f,g\rangle_\R=0\quad \text{for all}\,\,\, 
g\in H^\infty_+(\R/\langle\beta\rangle).
\label{eq-annih3}
\end{equation}
We recall that $g\in H^\infty_+(\R/\langle\beta\rangle)$ means that 
$g(x)=h(-\beta/x)$, for some function $h\in H^\infty_+(\R/2\Z)$. 
By the change-of-variables formula, we have
\begin{equation}
\langle f,g\rangle_\R=\int_\R f(x)g(x)\diff x=\int_\R f(x)h
\bigg(-\frac{\beta}{x}\bigg)\diff x=\beta
\int_\R f\bigg(-\frac{\beta}{x}\bigg)h(x)\frac{\diff x}{x^2}=\langle
\Jop_\beta f,h\rangle_\R,
\label{eq-annih4}
\end{equation}
where $\Jop_\beta:L^1(\R)\to L^1(\R)$ denotes the isometric transformation
\[
\Jop_\beta f(x):=\frac{\beta}{x^2}f\bigg(-\frac{\beta}{x}\bigg).
\]
From \eqref{eq-annih4} we see that \eqref{eq-annih3} is equivalent to 
having
\[
\Qop_2\Jop_\beta f\in H^1_{+,0}(\R/2\Z).
\]
It is easy to check that 
\[
f\in H^1(\R)\quad\Longrightarrow\quad \Qop_2 f,\,\Qop_2\Jop_\beta f\in 
H^1_{+,0}(\R/2\Z).
\]
Problem \ref{prob-1} asks whether, for $0<\beta\le1$, the reverse implication
holds: Is it true that
\begin{equation}
f\in L^1(\R)\,\,\,\text{and}\,\,\,\Qop_2 f,\,\Qop_2\Jop_\beta f\in 
H^1_{+,0}(\R/2\Z)\quad\Longrightarrow\quad f\in H^1_+(\R)?
\label{eq-impli1}
\end{equation}
We note that if we ask that, in addition, $f\in H^1_{\mathrm{real}}(\R)$, the 
implication holds, by the preceding argument.
\medskip

\noindent{\bf A related open problem.} It may shed light on \eqref{eq-impli1}
to formulate the analogous statement in the setting of $L^p(\R)$, for $0<p<1$.
From the well-known (quasi-triangle) inequality
\[
|z_1+\cdots+z_n|^p\le |z_1|^p+\cdots+|z_n|^p,\qquad 0<p\le 1,
\]
we quickly see that $\Qop_2:L^p(\R)\to L^p(\R/2\Z)$ is bounded.    
It remains to define $\Jop_{\beta,p}$. We put
\[
\Jop_{\beta,p}[f](x):=\beta^{1/p}|x|^{-2/p}\theta_p(x)\,
f\bigg(-\frac{\beta}{x}\bigg),
\]
where the phase factor $\theta_p(x)$ is defined as follows: $\theta_p(x):=1$ 
for $x>0$, and $\theta_p(x):=\e^{-\imag 2\pi/p}$ for $x<0$. 
It is well understood how one defines the Hardy spaces $H^p_+(\R)$ and 
$H^p_{+,0}(\R/2\Z)$ as closed subspaces of $L^p(\R)$ and $L^p(\R/2\Z)$, 
respectively, also for $0<p<1$. 
We are ready to formulate the general problem.

\begin{prob}
$(0<p\le1)$ For which positive $\beta$ is it true that
\begin{equation*}
f\in L^p(\R)\,\,\,\text{and}\,\,\,\Qop_2 f,\,\Qop_2\Jop_{\beta,p} f\in 
H^p_{+,0}(\R/2\Z)\quad\Longrightarrow\quad f\in H^p_+(\R)?
\end{equation*}
\label{prob-2}
\end{prob}
\medskip

\noindent{\bf Distorted lattice-crosses.}
We consider the set $\Lambda_{\alpha}^{\langle\xi^0\rangle}$ of Theorem 
\ref{thm-3.8}. and assume $\alpha,\beta,m$ are all positive with 
$\alpha\beta m^2\le 4\pi^2$.

\begin{proof}[Proof of Theorem \ref{thm-3.8}]
Let $\mu\in\mathrm{ACH}(\Gamma_m)$ have $\widehat\mu=0$ on 
$\Lambda_{\alpha}^{\langle\xi^0\rangle}$. By Theorem \ref{thm-3}, 
we have that $\widehat\mu$ vanishes on the set
\[
\bar\R^2_{--}\cup(\bar\R^2_{++}+\{\xi^0\}).
\]
In terms of the compressed measure $\bpi_1\mu\in\mathrm{M}(\R^\times)$, 
this is equivalent to having $\diff\bpi_1\mu(t)=f(t)\diff t$, where
$f\in H^1_-(\R)$ and $f/U_{\xi^0}\in H^1_+(\R)$, where $U_{\xi^0}$ is the
unimodular function
\[
U_{\xi^0}(t):=\e^{-\imag\pi[\xi_1^0t-m^2\xi_2^0/(4\pi^2 t)]},\qquad t\in\R.
\]
The given information allows us to conclude (e.g., we can use Morera's 
theorem) that $f$ has a holomorphic extension to $\C^\times=\C\setminus\{0\}$.
 
Now, if $\xi_1^0<0$, then the extension must decay too quickly as we approach
infinity in the upper half plane, so $f=0$ is the only possibility. 
If $\xi_1^0=0$, then still the point at infinity must be a removable 
singularity. If we look at the origin instead of infinity, we find that
if $\xi_2^0<0$, then the decay prescribed is too strong unless $f=0$. 
Moreover, if $\xi_2^0=0$, we get at least a removable singularity. 
So, if $\xi^0=(0,0)$, we get a removable singularity at the origin and at 
infinity, so by Liouville's theorem, $f$ must be constant, and the constant
is $0$, as $f\in H^1_-(\R)$. 
Nest, if $\xi_1^0>0$ and $\xi_2^0\ge0$, we may pick a non-trivial $f$ from a
Paley-Wiener space of entire functions (this is a closed subspace of $L^1(\R)$
of entire functions with the following properties: the functions are bounded 
in the lower half-plane, and have at most a given exponential growth in the 
upper half-plane). By applying the inversion $x\mapsto-1/x$, we can find
analogously non-trivial $f$ if $\xi_2^0>0$ and $\xi_1^0\ge0$.
The proof is complete.
\end{proof}

\section{Fourier uniqueness for a single branch of the hyperbola}

\noindent{\bf Dual formulation of the theorem.} We now turn to Theorem
\ref{thm-onebranch}, and observe that a scaling argument allows us to suppose
that
\[
\alpha=2,\quad m=2\pi.
\]
The dual formulation of Theorem \ref{thm-onebranch} now reads as follows. 
{\em The restriction to $\R_+$ of the functions 
\[
\e^{\imag2\pi jt},\,\,\,\e^{\imag\pi\beta k/t},\qquad j,k\in\Z,
\]
span a weak-star dense dense subspace of $L^\infty(\R_+)$ if and only if
$\beta<2$. Moreover, for $\beta=2$, the weak-star closure of the linear span
has codimension $1$ in $L^\infty(\R_+)$.}

\begin{proof}[Proof of Theorem \ref{thm-onebranch}]
Let $\nu\in\mathrm{AC}(\R_+)$; then $\nu$ may be written as 
$\diff\nu(t)=f(t)\diff t$, where $f\in L^1(\R_+)$. When needed, we think of 
$\nu$ and $f$ as defined to vanish on $\bar\R_-$. We suppose that
\begin{equation}
\int_0^{+\infty}\e^{\imag 2\pi jt}\diff\nu(t)=\int_0^{+\infty}
\e^{\imag2\pi\gamma k/t}\diff\nu(t)=0,\qquad j,k\in\Z,
\label{eq-fusb1}
\end{equation}
where $\gamma:=\beta/2$. We shall analyze the dimension of the space of 
solutions $\nu$, depending on the positive real parameter $\gamma$. 
We rewrite \eqref{eq-fusb1} in the form
\begin{equation}
\int_0^{+\infty}\e^{\imag2\pi jt}\diff\nu(t)=\int_0^{+\infty}
\e^{\imag2\pi kt}\diff\nu(\gamma/t)=0,\qquad j,k\in\Z,
\label{eq-fusb2}
\end{equation}
which we easily see is equivalent to having (cf. \cite{HM})
\[
\sum_{j\in\Z}\diff\nu(t+j)=\sum_{j\in\Z}\diff\nu
\bigg(\frac{\gamma}{t+j}\bigg)=0,\qquad t\in\R.
\]
Both expressions are $1$-periodic, so it is enough to require equality on
$[0,1[$ (we remove terms that are $0$):
\begin{equation}
\sum_{j=0}^{+\infty}\diff\nu(t+j)=\sum_{j=0}^{+\infty}\diff\nu
\bigg(\frac{\gamma}{t+j}\bigg)=0,\qquad t\in[0,1[.
\label{eq-fusb2.5}
\end{equation}
We single out the term with $j=0$, and obtain that
\begin{equation}
\diff\nu(t)=-\sum_{j=1}^{+\infty}\diff\nu(t+j),\qquad t\in[0,1[,
\label{eq-fusb3}
\end{equation}
and
\begin{equation}
\diff\nu(t)=
-\sum_{j=1}^{+\infty}\diff\nu\bigg(\frac{\gamma t}{\gamma+jt}
\bigg),\qquad t\in]\gamma,+\infty[.
\label{eq-fusb4}
\end{equation}
If we take absolute values, apply the triangle inequality, and integrate,
we get rather trivially from \eqref{eq-fusb3} that
\begin{equation}
\int_{[0,1[}\diff|\nu|(t)\le\sum_{j=1}^{+\infty}
\int_{[0,1[}\diff|\nu|(t+j)=\int_{[1,+\infty[}\diff|\nu|(t),
\label{eq-fusb5}
\end{equation}
and from \eqref{eq-fusb4} that
\begin{equation}
\int_{[\gamma,+\infty[}\diff|\nu|(t)\le
\sum_{j=1}^{+\infty}\int_{[\gamma,+\infty[}
\diff|\nu|\bigg(\frac{\gamma t}{\gamma+jt}
\bigg)=\sum_{j=1}^{+\infty}\int_{[\frac{\gamma}{j+1},\frac{\gamma}{j}[}
\diff|\nu|(t)=\int_{]0,\gamma[}\diff|\nu|(t).
\label{eq-fusb6}
\end{equation}
For $0<\gamma\le1$, we may combine \eqref{eq-fusb5} and \eqref{eq-fusb6}, 
to arrive at
\begin{equation}
\int_{[0,1[}\diff|\nu|(t)\le\int_{[1,+\infty[}\diff|\nu|(t)\le
\int_{[\gamma,+\infty[}\diff|\nu|(t)\le\int_{]0,\gamma]}\diff|\nu|(t),
\label{eq-fusb7}
\end{equation}
which is only possible if we have equality everywhere in \eqref{eq-fusb7}.
But then $|\nu|$ takes no mass on the interval $[\gamma,1]$, and we must also
have $(0<\gamma\le1)$
\begin{equation}
\diff|\nu|(t)=\sum_{j=1}^{+\infty}\diff|\nu|(t+j),\qquad t\in[0,1[,
\label{eq-fusb8}
\end{equation}
and
\begin{equation}
\diff|\nu|(t)=
\sum_{j=1}^{+\infty}\diff|\nu|\bigg(\frac{\gamma t}{\gamma+jt}
\bigg),\qquad t\in]\gamma,+\infty[.
\label{eq-fusb9}
\end{equation}
Moreover, for some constant $\zeta\in\C$ with $|\zeta|=1$,
we must also have $(0<\gamma\le1)$
\begin{equation}
\diff\nu(t)=\zeta\diff|\nu|(t), \qquad t\in[0,1[,
\label{eq-fusb10}
\end{equation}
and 
\begin{equation}
\diff\nu(t)=-\zeta\diff|\nu|(t), \qquad t\in[1,+\infty[.
\label{eq-fusb11}
\end{equation}
For $0<\gamma\le1$, this allows us to focus on the positive measure 
$\diff|\nu|$. Again for $0<\gamma\le1$, we may combine \eqref{eq-fusb8} and 
\eqref{eq-fusb9} and obtain as a result that
\begin{equation}
\diff|\nu|(t)=
\sum_{j,k=1}^{+\infty}\diff|\nu|\bigg(\frac{\gamma(t+j)}{\gamma+k(t+j)}
\bigg),\qquad t\in]0,1[.
\label{eq-fusb12}
\end{equation}
For $x\in\R$, let $\{x\}_1$ be the fractional part of $x$; more precisely,
$\{x\}_1$ is the number in $[0,1[$ such that $x-\{x\}_1\in\Z$. 
We define $U_\gamma:[0,1[\to[0,1[$ as follows:
$U_\gamma(0):=0$, and 
\begin{equation}
U_\gamma(x):=\{\gamma/x\}_1,\qquad x\in]0,1[.
\label{eq-fusb13}
\end{equation}
As we already observed, $|\nu|$ takes no mass on $[\gamma,1]$. If we integrate
the left hand side of \eqref{eq-fusb12} on $[\gamma,1]$ to get $0$, we should
obtain $0$ from the right hand side as well. But integration of the right
hand side on $[\gamma,1]$ computes the $|\nu|$-mass of the set
\[
E_\gamma(2):=\{t\in[0,1[:\,\,U_\gamma^{2}(t)\in[\gamma,1]\},
\]
where $U_\gamma^2=U_\gamma\circ U_\gamma$, the composition square.
So $|\nu|$ takes no mass on $E_\gamma(2)$. By iterating this argument, we
see that $\mu$ assumes no mass on all sets of the form
\[
E_\gamma(2n):=\{t\in[0,1[:\,\,U_\gamma^{2n}(t)\in[\gamma,1]\},
\qquad n=1,2,3,\ldots.
\]
If $0<\gamma<1$, the union of all the sets $E_\gamma(n)$, $n=1,2,3,\ldots$,
has full Lebesgue mass, which no place for the mass of $|\nu|$, and we get that
$|\nu|([0,1[)=0$. By \eqref{eq-fusb7}, we get that $|\nu|(\R)=0$, that is,
$\nu=0$ identically. 

The case $\gamma=1$ is a little different. Then \eqref{eq-fusb12} asserts that
$|\nu|$ is an invariant measure for $U_1^2$, the square of the standard 
Gauss map \cite{CFS}.
As $U_1$ is ergodic with respect to the absolutely continuous probability 
measure 
\[
\diff\varpi(t):=\frac{\diff t}{(1+t)\log2}, 
\]
we conclude that $|\nu|$ must be of the form 
\[
\diff|\nu|(t)=C_1\diff\varpi(t),\qquad t\in[0,1[,
\]
for some real constant $C_1\ge0$. The analogous argument based on the interval
$[1,+\infty[$ in place of $[0,1[$ gives that
\[
\diff|\nu|(t)=C_2\diff\varpi(1/t)=\frac{C_2\diff t}{t(1+t)\log2},
\qquad t\in[1,+\infty[.
\]
We obtain that $\diff\nu$ must be a complex constant multiple of the measure
\[
1_{[0,1[}(t)\frac{\diff t}{1+t}-1_{[1,+\infty[}(t)\frac{\diff t}{t(1+t)}.
\]
This measure meets \eqref{eq-fusb3} and \eqref{eq-fusb4}, so we really have
a one-dimensional annihilator for $\gamma=1$.

Finally, we need to consider $\gamma>1$, and supply a non-trivial $\nu\in
\mathrm{AC}(\R_+)$ with \eqref{eq-fusb2.5}. In this case,
the Gauss-type map $U_\gamma$ given in \eqref{eq-fusb13} is uniformly 
expanding, and therefore, it has a non-trivial absolutely continuous 
invariant probability measure on $[0,1]$, which we call $\varpi_\gamma$ 
(cf. \cite{CFS}, p. 169, and \cite{ChR}, \cite{ChD}). 
We extend $\varpi_\gamma$ to $\R_+$ trivially by putting it equal to the zero 
measure on $\R_+\setminus[0,1]=]1,+\infty[$. Being invariant, 
$\varpi_\gamma$ has the property
\[
\diff\varpi_\gamma(t)=\sum_{j=1}^{+\infty}
\diff\varpi_\gamma\bigg(\frac{\gamma}{t+j}\bigg),\qquad t\in[0,1].
\]
We put
\[
\diff\nu(t):=\diff\varpi_\gamma(t)-\diff\varpi_\gamma(\gamma/t),\qquad
t\in\R_+,
\]
so that $\nu$ gets to have the symmetry property
\[
\diff\nu(t)=-\diff\nu(\gamma/t),\qquad t\in\R_+.
\]
It is now a simple exercise to verify that $\nu$ meets \eqref{eq-fusb2.5},
which completes the proof.
\end{proof}

\begin{rem}
\label{rem-g1}
In case $\gamma>1$, it is of interest to know how to construct more general 
measures $\nu\in\mathrm{AC}(\R_+)$ with \eqref{eq-fusb2.5}. We could try with
$\nu$ of the form 
\[
\diff\nu(t)=\diff\omega_1(t)-\diff\omega_1(\gamma/t)-\diff\omega_2(t),
\]
where $\omega_1$ is supported on $[0,1]$, while $\omega_2$ is supported on 
$[1,\gamma]$. We require that, in addition, 
$\diff\omega_2(\gamma/t)=-\diff\omega_2(t)$. Then $\nu$ has the symmetry 
property $\diff\nu(\gamma/t)=-\diff\nu(t)$, and we just need to check whether
\[
\sum_{j=0}^{+\infty}\diff\nu(t+j)=0,\qquad t\in]0,1[.
\]
We obtain the equation
\[
\diff\omega_1(t)=\sum_{j=0}^{+\infty}\diff\omega_1
\bigg(\frac{\gamma}{t+j}\bigg)+\sum_{j=1}^{]\gamma[}\diff\omega_2(t+j),
\qquad t\in]0,1[
\]
where $]\gamma[$ denotes the largest integer $<\gamma$. In particular, if
$1<\gamma\le2$, this equation reads
\[
\diff\omega_1(t)=\sum_{j=0}^{+\infty}\diff\omega_1
\bigg(\frac{\gamma}{t+j}\bigg)+\diff\omega_2(t+1).
\qquad t\in]0,1[
\]
This equation is perturbation of the invariant measure equation (which is
obtained for $\omega_2=0$), and one would expect that there should exist 
many solutions $\omega_1,\omega_2$ (cf. \cite{CHM}).
\end{rem}
\medskip

\noindent{\bf Analysis of the critical case $\alpha\beta m^2=16\pi^2$.}
Without loss of generality, we may take
\[
\alpha=\beta=2,\,\,\, m=2\pi,
\]
which corresponds to $\gamma=1$ in the above proof of Theorem 
\ref{thm-onebranch}. We now look at the cause of the defect $1$, the 
one-dimensional subspace spanned by the measure
\[
\diff\nu(t)=
1_{[0,1[}(t)\frac{\diff t}{1+t}-1_{[1,+\infty[}(t)\frac{\diff t}{t(1+t)},
\]
as we see from the proof of Theorem \ref{thm-onebranch}. This measure has 
the symmetry property $\diff\nu(1/t)=-\diff\nu(t)$, which means that
\[
\int_{\R_+}\e^{\imag 2\pi x/t}\diff\nu(t)=-\int_{\R_+}
\e^{\imag 2\pi x t}\diff\nu(t),\qquad x\in\R. 
\]
We will need to compute the one-dimensional Fourier transform 
\[
\widehat\nu(x):=\int_{\R_+}\e^{\imag 2\pi x t}\diff\nu(t),\qquad x\in\R.
\]
We quickly find that
\[
\widehat\nu(x):=(1-\e^{\imag 2\pi x})\int_{0}^{+\infty}
\e^{\imag 2\pi x t}\frac{\diff t}{t+1},\qquad x\in\R,
\]
where the integral on the right hand side is understood in the generalized 
Riemann sense. 

\begin{proof}[Proof of Corollary \ref{cor-onebranch}]
By symmetry, we may take $\xi^0\in\R\times\{0\}$. It will be enough to
establish that 
\[
\widehat\nu(x)\ne0,\qquad x\in\R\setminus\Z.
\] 
It will be sufficient to obtain that
\[
\int_{1}^{+\infty}\e^{\imag 2\pi x t}\frac{\diff t}{t}
=\int_{1}^{+\infty}\cos(2\pi x t)\frac{\diff t}{t}+
\imag\int_{1}^{+\infty}\sin(2\pi x t)\frac{\diff t}{t}\ne0,
\qquad x\in\R^\times.
\]
The real part of this expression equals
\[
\int_{1}^{+\infty}\cos(2\pi x t)\frac{\diff t}{t}
=\int_{|x|}^{+\infty}\frac{\cos y}{y}\,\diff y=-\ci(|x|),
\]
whereas the imaginary part equals
\[
\int_{1}^{+\infty}\sin(2\pi x t)\frac{\diff t}{t}=
\sign(x)\int_{|x|}^{+\infty}\frac{\sin y}{y}\,\diff y=-\sign(x)\si(|x|);
\]
the $\sign$ function was defined in Section \ref{sec-Hilb}, and the integral
expression can be thought of as defining the rather standard functions 
``si'' and ``ci''.  It is well-known that the parametrization
\[
\ci(x)+\imag\si(x),\qquad 0<x<+\infty,
\]
forms the {\em Nielsen} (or sici) spiral which converges to the origin as 
$x\to+\infty$, and whose curvature is proportional to $x$ 
(see, e.g. \cite{PlM}). 
In particular, the spiral never intersects the origin, which does it. 
\end{proof}

\section{Open problems in higher dimensions}

\noindent{\bf The Klein-Gordon equation in dimension $d$.} In space
dimension $d>1$, we consider a solution $u$ to \eqref{eq-KG100} of the form
\[
u(t,x)=\hat\mu(t,x):=\int_{\R^{d+1}}\e^{\pi\imag(\tau t+\langle x,\xi\rangle)}
\diff\mu(\tau,\xi),
\]
where $\mu$ is a complex-valued finite Borel measure, and $t,\tau\in\R$,
  $x,\xi\in\R^d$, and 
\[
\langle x,\xi\rangle=x_1\xi_1+\cdots+x_d\xi_d.
\]
The assumption that $u$ solves the Klein-Gordon equation means that
\[
\bigg(\tau^2-|\xi|^2-\frac{\EP}{\pi^2}\bigg)\diff\mu(\tau,\xi)=0
\]
as a measure on $\R^{d+1}$, which we see is the same as having
\begin{equation*}
\supp\mu\subset\Gamma_{m}(d):=\bigg\{(\tau,\xi)\in\R\times\R^d:
\,\tau^2-|\xi|^2=\frac{\EP}{\pi^2}\bigg\}.
\end{equation*}
The set $\Gamma_m(d)$ is a two-sheeted $d$-dimensional hyperboloid. Let
\[
\Gamma_m^+(d):=\{(\tau,\xi)\in\Gamma_m:
\,\tau>0\},\qquad 
\Gamma_m^-(d):=\{(\tau,\xi)\in\Gamma_m:
\,\tau<0\},
\]
be the two connectivity sheets of the hyperboloid $\Gamma_m(d)$. We equip
$\Gamma_m$ with $d$-dimensional surface measure, and require of $\mu$ that it
be absolutely continuous with respect to this surface measure.
\medskip

\noindent{\bf Light cones.} 
We consider the light cone emanating from the origin: 
\[
{\mathrm Y}_0:=\big\{(t,x)\in \R\times\R^{d}:\,\,|x|=|t|\big\}.
\]
The light cone is a characteristic surface for the Klein-Gordon equation.
For any $\varepsilon\ge0$, the surface
\[
{\mathrm Y}_0(\varepsilon):=\big\{(t,x)\in \R\times\R^{d}:\,\,|x|=|t|
+\varepsilon\big\}
\]
is characteristic as well. In connection with their study of the event 
horizon of Kerr black holes, Ionescu and Klainerman \cite{IK} showed 
(for $\varepsilon>0$) that if the function $u$ -- which solves the 
Klein-Gordon equation -- vanishes on ${\mathrm Y}_0(\varepsilon)$, then 
$u=0$ for all $(t,x)$ with $|x|\ge|t|+\varepsilon$ (so we get suppression in 
the space-like direction); compare also with \cite{KRS} and \cite{Nir}.
Klainerman (private communication) has indicated that this should be true
for $\varepsilon=0$ as well. But then we should expect $\mathrm{Y}_0$ to be
a uniqueness set for $u$, as there is no width to the waist of $\mathrm{Y}_0$
which could be the source for a wave. So, we suppose for the moment that it
has been established that $\mathrm{Y}_0$ is a uniqueness set for $u$. 
Then it makes sense to ask for (small) subsets of $\mathrm{Y}_0$ that are
sets of uniqueness, too. This is what Theorem \ref{thm-1} supplies in $d=1$. 
In analogy with Theorem \ref{thm-onebranch}, we would ask for even smaller
subsets of $\mathrm{Y}_0$ that are sets of uniqueness for $u$, provided that
the Borel measure $\mu$ (which $u$ is the Fourier transform of) is supported
on the branch $\Gamma_m^+(d)$.


\end{document}